\documentclass{article}

\usepackage{latexsym}
\usepackage{amssymb,amsmath,amsopn}
\usepackage[dvips]{graphicx}
\usepackage{color,epsfig}
\usepackage{url}
\usepackage{multicol}
\usepackage{amsthm}

\pdfpagewidth \paperwidth
\pdfpageheight \paperheight
\newtheorem*{remark}{Remark}
\newtheorem{theorem}{Theorem}
\newcommand*{\twoheadrightarrowtail}{\mathrel{\rightarrowtail\kern-1.9ex\twoheadrightarrow}}

\begin{document}
\title{\textbf{The generalized \\Sierpi\'{n}ski Arrowhead Curve}}

\author{Kaszanyitzky, Andr\'as \\ kaszi75@gmail.com}
\date{}
\maketitle

\textbf{Abstract.} In \emph{Section 1} we define and observe a special \emph{Hamiltonian-path} and a special tiling-path on a checked triangular grid which is related to the \emph{generalized Sierpi\'{n}ski Gasket}. We describe these paths with strings formed by the absolute direction codes of their edges and prove their existence on triangular grids in any order of $n$ by a constructing algorithm. We find and prove a bijective relation between them with a transformation table, which ensures their unambiguous transformability into each other.

In \emph{Sections 2} and \emph{3} we extend our bijection to larger graphs and use these bijective pairs as generator curves of the \emph{generalized Sierpi\'{n}ski Arrowhead Curve}. The cardinality of these generator curves specifies a new integer sequence. Our curves keep their basic properties (simplicity, self-avoidance) independently of the transformation between node-rewriting and edge-rewriting methods unlike other recursive curves. They continuously map the fractal family in many ways and differ from \emph{Pascal Triangle modulo $n$} patterns in compound orders of $n$. We produce these symmetric curves by asymmetric paths. We also find their \emph{Hausdorff dimension}. Triangular numbers and their powers are also mentioned.

In \emph{Section 4} we also discuss how to transform our paths into the \emph{Lindenmayer-system} to draw the recursive curves in an easy way.

\section{The bijective relation between paths and tilings}

In this section we define two kinds of special paths related to the same generator pattern on the triangular grid. We observe their transformability into each other and prove their bijective relation.

\subsection{Checked generator patterns}

Let us consider an equilateral triangle with sides divided into $n$ equal pieces. Connecting the dividing points with line-segments parallel to the sides, we get a natural partitioning of the $n^2$ smaller coincident subtriangles by colouring the tiles that face upwards (like the original triangle) dark, and colouring the rest of the subtriangles white.

By infinitely replacing all the dark tiles with contracted copies of the checked generator pattern, we get a symmetric fractal, the \emph{generalized Sierpi\'{n}ski Gasket}.  This 2-dimensional checked generator pattern of order $n$ in the $k$-th approximation (repeating this transformation $k$ times) usually denoted by $SG_{2,n}(k)$. We will be referring to the generator pattern in a shorter form by $F_n$, and the approximations of the fractal by $F_n(k)$. See \emph{Figure 1}.

\begin{figure}[ht]
\centering
\includegraphics[width=0.9\textwidth]{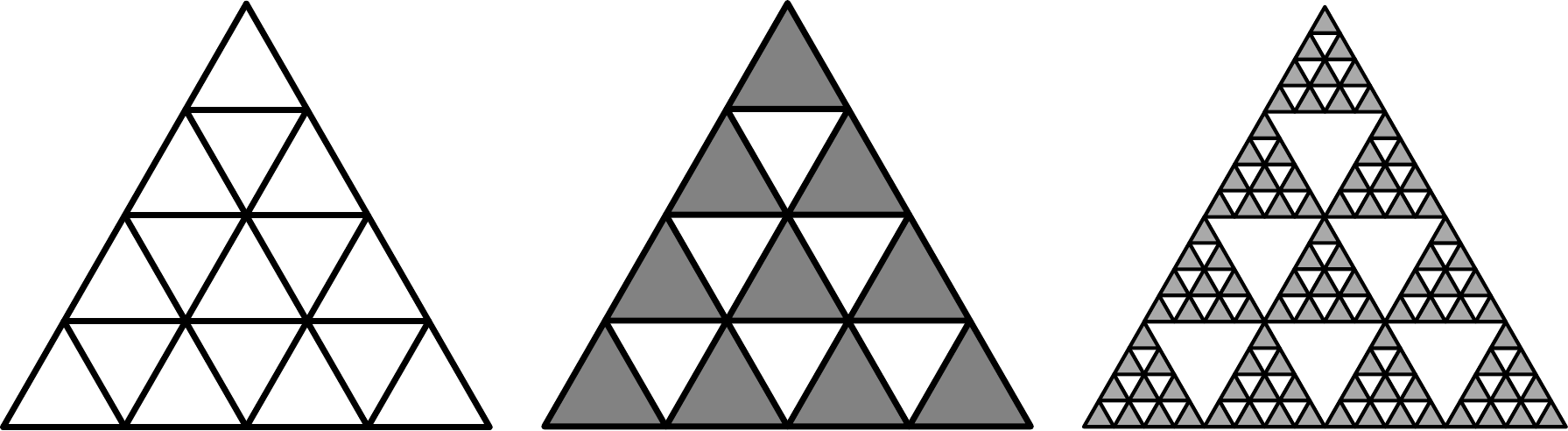}
\centerline{}
\centerline{\bf Figure 1. Triangular grid of order 4 with 16 tiles,}
\centerline{\bf $F_4$ checked generator pattern and the 2nd approximation: $F_4(2)$.}
\end{figure}

This fractal family keeps the original properties of the pattern for any order of $n$, unlike \emph{Pascal Triangle modulo $n$} patterns which keep these properties only in prime orders and contain other patterns in the white subtriangles which face downwards.

Our generator pattern $F_n$ contains $n^2$ subtriangles, $T_{n-1}$ white tiles, $T_n$ dark tiles and $T_{n+1}$ grid points as consecutive triangular numbers. The centroids of the monochromatic tiles also form a triangular grid. We will only be using the dark tiles. Their centroids form the \emph{inscribed grid}, and their corners form the \emph{overall grid}. See \emph{Figure 2}.

\begin{figure}[ht]
\centering
\includegraphics[width=0.7\textwidth]{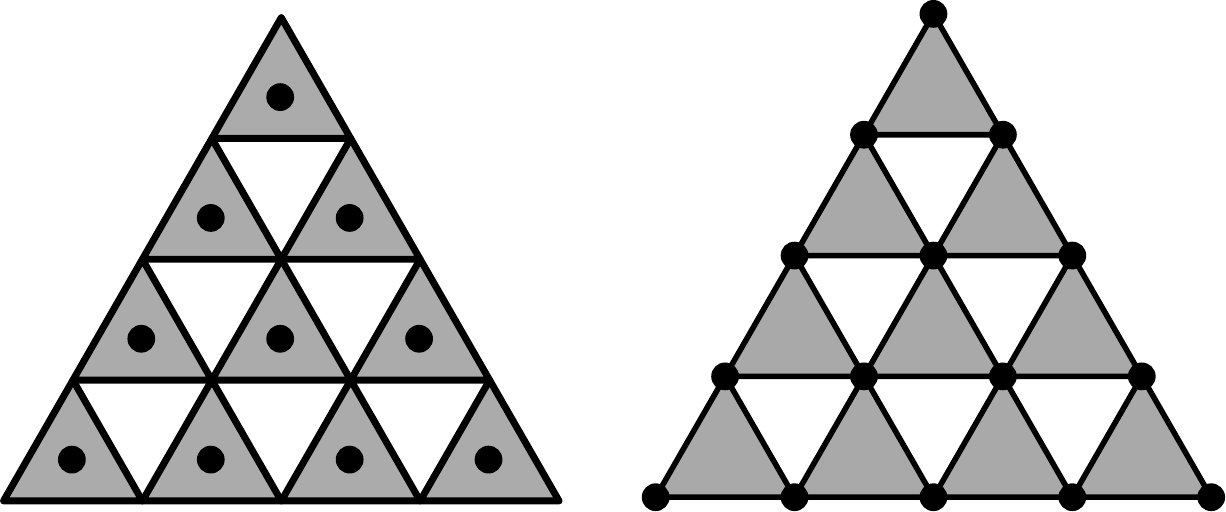}
\centerline{}
\centerline{\bf Figure 2. Generator pattern of order 4 consists of $T_4$ dark tiles.}
\centerline{\bf Dots define the inscribed grid (left) and the overall grid (right).}
\centerline{\bf They consist of $T_4=10$ and $T_5=15$ grid points.}
\end{figure}

\subsection{Paths and permutations of the dark tiles}

By connecting the $T_n$ grid-points of the inscribed grid with the shortest possible edges in all the self-avoiding ways we get Hamiltonian-paths, \emph{H-paths} denoted by $H_n$.

Let us consider a self-avoiding tiling-path on the overall grid called \emph{S-path} (referring to Sierpi\'{n}ski), denoted by $S_n$ which consists of $T_n$ edges and in which all the edges must be lying on different dark subtriangles. For practical reasons we will be using the notation of McKenna: marking the tiles with little ticks in the middle of the edges [McK94]. See the right side of \emph{Figure 3}.

\begin{figure}[ht]
\centering
\includegraphics[width=0.7\textwidth]{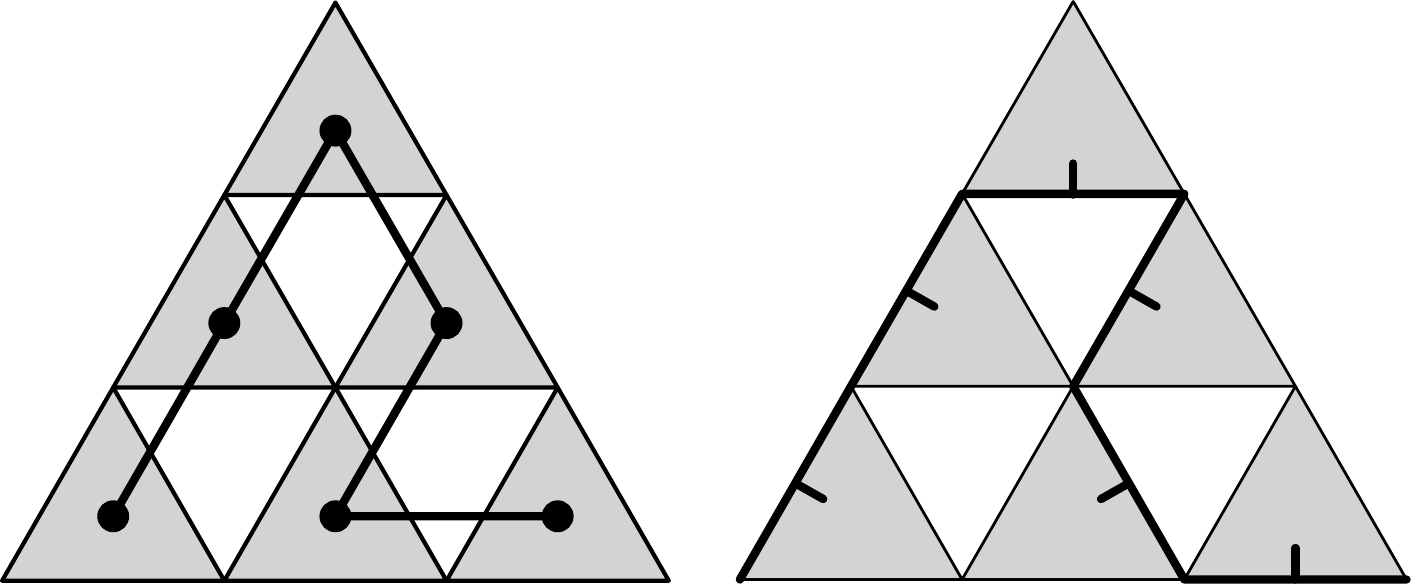}
\centerline{}
\centerline{\bf Figure 3. Example of an $H_3$-path and the corresponding $S_3$-path}
\centerline{\bf with the same permutation of the dark tiles.}
\end{figure}

Both paths have directed edges which are originating from the leftmost grid point and terminating in the rightmost grid point.

H-paths and S-paths are both connect node-neighbour dark subtriangles and describe permutations of the dark tiles, but these paths have different cardinality.

\subsection{Absolute direction code of the edges}

We will describe both paths with strings of the absolute direction codes of their edges.

Let us define a triangular grid $T_n$ with $A,B,C$ corners where $AB$ is a polar axis. The third corner $C$ is located in positive direction (counterclockwise) on our polar coordinate system.

We are searching for paths from $A$, the originating leftmost grid point (pole) to $B$, the terminating rightmost grid point (the end point of the polar axis).

Consider the edges of the path translated into the pole. The only unique property is their direction. Let us denote this absolute direction by $d$, where ${0}\leq{d}\leq{5}$ and a direction code means $d\dfrac{\pi}{3}$ radian. See \emph{Figure 4}.

\begin{figure}[ht]
\centering
\includegraphics[width=0.7\textwidth]{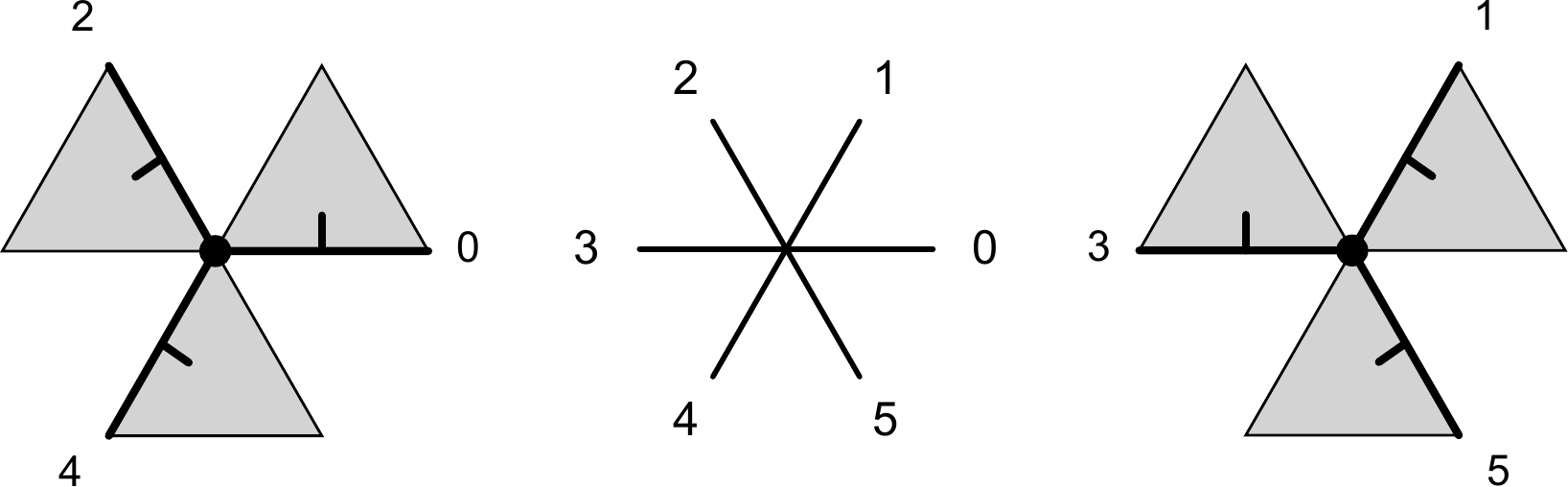}
\centerline{}
\centerline{\bf Figure 4. Absolute direction codes of H-paths (middle)}
\centerline{\bf and S-paths (left and right sides) in which}
\centerline{\bf parity specifies the facing of the tiles.}
\end{figure}

We can describe the consecutive edges of our paths with a string which consists of these direction codes in both $H$-paths and $S$-paths.

It is obvious that there are no turn backs in paths, so for every \linebreak $d_i$ and $d_{i+1}$, $\mid{d_{i+1}-d_{i}} \mid \neq 3$ stands for both paths.

The definition of the S-path also excludes turn backs to the same dark tile, therefore consecutive edges $(d_i, d_{i+1})$ of the S-path cannot be the following: 

\begin{equation*}
d_{i+1}\not\equiv
\begin{cases}
(d_i +2) \mod 6 & \qquad \text{if $d_i$ is even} \\
(d_i +4) \mod 6 & \qquad \text{if $d_i$ is odd}
\end{cases}
\end{equation*}

\subsection{Transformability of the paths into each other}

Here we introduce and observe the transformations between the paths in both directions. We will prove that there is a bijection between a subset of $H$-paths and $S$-paths which are related to the same generator pattern.

The transformation from an S-path into a Hamiltonian-path means we connect the centroids of all the dark tiles exactly in the same order as the edges of the S-path touched them. This transformation is always possible and results exactly one H-path which is different for every S-path. This is an injective relation.

The transformation from an H-path into an S-path means we connect all the dark tiles from the leftmost point of the overall grid through the contact points of the node-neighbour dark tiles to the rightmost point with tiling-edges exactly in the same order as the H-path connected them. Sometimes it is not possible, but any of the S-paths have a preimage in the set of H-paths, therefore it is a surjective relation.

\begin{remark}
There are more Hamiltonian-paths on the inscribed grid than S-paths on the overall grid: $\mid{H_n} \mid \geqslant \mid{S_n} \mid$, because some of them are impossible to transform into an S-path.
\end{remark}

Consider 3 dark node-neighbour tiles which share on one common node in the middle. By connecting their centroids we can always form a Hamiltonian-path, but it is impossible to do so with S-paths. We cannot draw 3 consecutive edges which touch 3 different tiles in this arrangement because the 3 tiles have only one contact point instead of two. See \emph{Figure 5} where we cannot continue the S-path with an edge which touches the third tile signed by an X.

\begin{figure}[ht]
\centering
\includegraphics[width=0.5\textwidth]{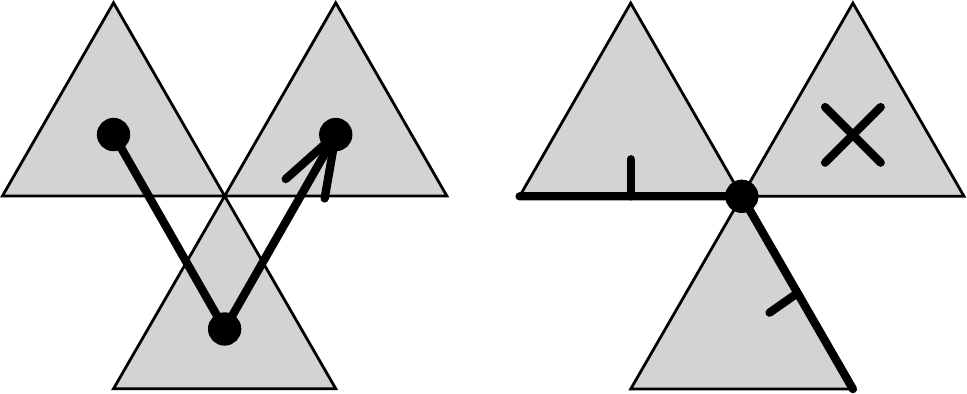}
\centerline{}
\centerline{\bf Figure 5. Impossibility of the transformation}
\centerline{\bf from an H-path into an S-path.}
\end{figure}

These injective $(S_n \rightarrowtail H_n)$ and surjective $(H_n \twoheadrightarrow S_n)$ relations mean there is a bijection between a subset of $H_n$ and $S_n$. We call this subset the Well-formed Hamiltonian-paths (W-paths), and referring to them by $W_n$.

\subsection{Special Hamiltonian paths}

\begin{theorem}
W-paths on the inscribed grid and S-paths on the corresponding overall grid are in bijective relation if they belongs to the same generator pattern: $W_n \twoheadrightarrowtail S_n$.
\end{theorem}

By excluding some undesirable turns we can define W-paths, a special subset of H-paths which can be bijectively transformed into S-paths and vice versa.

We use new notations ($s_n$, $h_n$ and $w_n$) for the strings of the direction codes referring to the paths of order $n$ denoted by the same capital letters.

Let us consider all possible Hamiltonian-paths ($H_n$) on the inscribed grid, with the absolute direction codes described above. We would like to observe the turns of the edge-pairs in the H-path, therefore we consider a basic direction ($0$) from the left to the right same as the position of the originating and the terminating point and we supplement our $h_n$ string with a leading and an ending zero.

The forbidden turns as $(d_i, d_{i+1})$ number pairs are the following: 

\begin{equation*}
d_{i+1}\not\equiv
\begin{cases}
(d_i +4) \mod 6 & \qquad \text{if $d_i$ is even} \\
(d_i +2) \mod 6 & \qquad \text{if $d_i$ is odd}
\end{cases}
\end{equation*}

where $d_i$ and $d_{i+1}$ are neighbouring elements of the direction code string $h_n$. Geometrically, this means that the next edge cannot turn $120^{\circ}$ to the right after an even direction and it cannot turn $120^{\circ}$ to the left after an odd direction. You can see a wrong turn on the left side of \emph{Figure 6} at the 2nd and the 3rd edges of the H-path, as direction code pairs 51. We call these kinds of strings consisting of only well-formed turns $w_n$ strings in that case if this condition is true.

\begin{figure}[ht]
\centering
\includegraphics[width=0.7\textwidth]{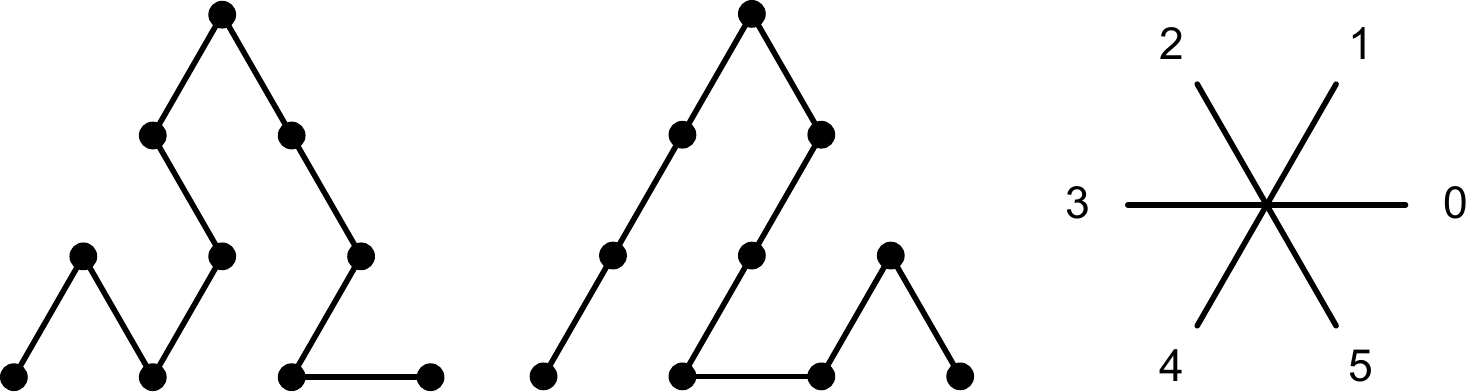}
\centerline{}
\centerline{\bf Figure 6. H-path ($h_4=151215540$) and W-path ($w_4=111544015$).}
\end{figure}

\subsection{Trivial W-paths}

\begin{remark}
Well-formed paths exist for all orders.
\end{remark}

We will be showing a constructing algorithm for a trivial W-path which proves their existence on the inscribed grid for all orders.

Let $j$ be equal to 1. 

The iterated steps are the following: head $(n-j)$ times up-right (1), then head down-right $1$ time (5). Increment $j$, then head $(n-j)$ times down-left (4), then go to the right $1$ time (0), and increment $j$.

If $n$ is odd, we have to iterate our steps $\dfrac{n-1}{2}$ times, otherwise $\dfrac{n-2}{2}$ times, and in this case we have to add 2 more edges (an arrowhead shape) to the string in the directions up-right (1) and down-right (5) after the iterations have finished.

For the $n=4$ case of our trivial W-path, see the middle image of \emph{Figure 6}.

When $n=6$, our trivial W-path as a string (using the absolute direction codes we previously described) will be the following:  

\smallskip

\centerline{$w_6=11111544440111544015$}

\smallskip

This proves that we can always construct at least one W-path on our inscribed grid for any order of $n$.

\subsection{Transformation table of the paths}

In W-paths, there are no untransformable turns. They have the same cardinality as S-paths: $\mid{W_n} \mid = \mid{S_n} \mid$. They describe the same permutations of the dark tiles. The connecting points of the dark tiles in an S-path can be transformed into the edges of the unique corresponding W-path and vice versa. We will be showing this property and prove the bijection by a transformation table. See \emph{Table 1, Figure 4} and \emph{Figure 7}.

We always have to supplement our $w_n$ string with a leading and an ending zero to begin transforming a W-path into an S-path: $0 w_n 0 \rightarrow s_n$. Then every pair $(a,b)$ of the supplemented $w_n$ string (first and second digits, second and third digits, etc.) give a new direction code, a digit of the $s_n$ string, by reading \emph{Table 1} by rows and columns $(a,b)$. These pairs describe two contact points between 3 dark tiles ($P_1 P_2, P_2 P_3,  P_3 P_4$ on \emph{Figure 7}). The middle tile has an entering and an exiting point which are in a one-to-one correspondence with the direction code of the actual edge in the S-path. (See the left side and the middle of \emph{Figure 7} and for the direction codes see \emph{Figure 4}.)

\[
\renewcommand{\arraystretch}{1.3}
\begin{tabular}{|c||c|c|c|c|c|c|}
\hline
\bf
$\blacksquare$ & {\bf $b=0$} & \bf {1} & \bf {2} & \bf {3} & \bf {4} & \bf {5} \\ \hline \hline
\bf {$a=0$} & 0 & 1 & 1 & $\blacksquare$ & $\blacksquare$ & 0 \\ \hline
\bf 1 & 0 & 1 & 1 & $\blacksquare$ & $\blacksquare$ & 0 \\ \hline
\bf 2 & $\blacksquare$ & 2 & 2 & 3 & 3 & $\blacksquare$ \\ \hline
\bf 3 & $\blacksquare$ & 2 & 2 & 3 & 3 & $\blacksquare$ \\ \hline
\bf 4 & 5 & $\blacksquare$ & $\blacksquare$ & 4 & 4 & 5 \\ \hline
\bf 5 & 5 & $\blacksquare$ & $\blacksquare$ & 4 & 4 & 5 \\ \hline
\end{tabular}
\]

\smallskip 

\centerline {\bf \emph{Table 1.} Transformation table for changing the direction codes} \centerline {\bf from a W-path into an S-path and vice versa.}

\begin{figure}[ht]
\centering
\includegraphics[width=0.9\textwidth]{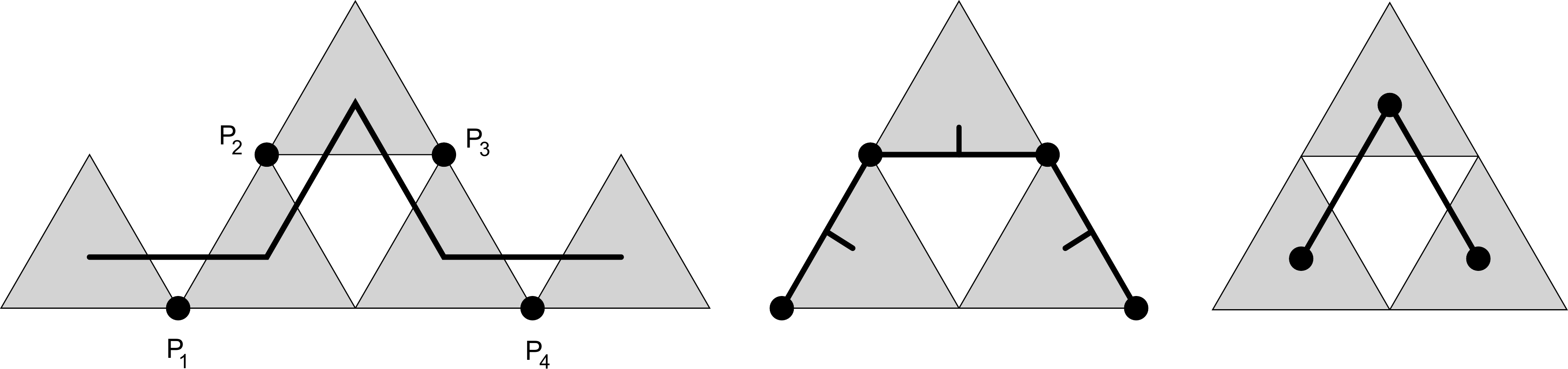}

\centerline{\bf Figure 7. Changing a supplemented W-path into an S-path}
\centerline{\bf by reading Table 1 as values defined by $(a,b)$}
\centerline{\bf and vice versa by reading it in reverse order $(b,a)$.}
{\bf $0 w_2 0=0150 \qquad \rightarrow \qquad s_2=105 \qquad \rightarrow \qquad w_2=15$.}
\end{figure}

To change the direction codes from an S-path to a W-path ($s_n \rightarrow w_n$) we have to read \emph{Table 1} by columns and rows $(b,a)$ because this is the transpose of our matrix.

In this case consecutive edge-pairs of the S-path $(b,a)$ gives us the edges of the inscribed grid (digits of the $w_n$ string). Consecutive digit-pairs of the $s_n$ string show the unique possible arrangement of 2 neighbour dark tiles connected by a contact point ($P_2, P_3$). We can connect their centroids only in one way to get the edges of the W-path. (See the middle and the right side of \emph{Figure 7}.)

The definitions of the paths with the rules and properties of their transformations, the presented illustrations, the constructing algorithm of the trivial paths, the transformation table and all these considerations prove \emph{Theorem 1}.

\section{Triangular self-avoiding recursive curves}

In this section we will show that these bijective pairs of the paths are the generator curves of the \emph{generalized Sierpi\'{n}ski Arrowhead Curve}. We describe their construction methods. We will call our fractal approximating extended grids the \emph{inscribed graph} and the \emph{overall graph}.

\subsection {FASS-curves}

In 1890, \emph{Peano} discovered the first \emph{space-filling curve}, which was the first continuous geometrical example for \emph{Cantor}'s surprising theorems.

In 1994, \emph{Douglas M. McKenna} proved that unlike on a square grid, simple triangular Peano-curves (FASS-curves) do not exist [McK94]. I thought non space-filling curves are worth observing, and I found a new family of simple self-avoiding recursive curves in 2009.

For further study of our subject, simple \emph{recursive curves} or \emph{Peano-curves}, also known as \emph{FASS-curves}, which have space-Filling, self-Avoiding, Simple and self-Similar properties, and the related \emph{Lindenmayer-systems} with their rewriting methods we recommend reading the followings: [PL90, PLF91, Sa94, M82].

\subsection{Rewriting methods}

Recursive curves are usually described by a formal language called L-system, named after the Hungarian theoretical biologist \emph{Aristid Lindenmayer}, who lived and worked in Utrecht, the Netherlands.

There are 2 different ways to make a recursive curve in the Lindenmayer-system: the node-rewriting (NR) and the edge-rewriting (ER) method.

It is an obvious but often overlooked fact that \emph{node-rewriting} FASS-curves are Hamiltonian-paths in any approximation. In an NR path we keep the edges and substitute \emph{the nodes} in all approximations with the transformed copy of the original path.

\emph{Edge-rewriting} curves work like tessellations. Edges of this self-avoiding path do not visit all the grid points, but they touch every little tile exactly once. We subsitute \emph{the edges} in all approximations with the transformed copy of the original path.

Both of our paths can be made extendable with a simple self-similar definition to larger graphs. A \emph{W-path} on the inscribed grid became to \emph{NR recursive curve} and results Hamiltonian-paths on the inscribed graph. An \emph{S-path} on the overall grid became to \emph{ER recursive curve} and results tiling-paths on the overall graph. They are not FASS-curves because they are \emph{non space-filling paths}, but they describe and fill the same fractal pattern $F_n(\infty)$ in an infinite approximation. 

\subsection{Approximately space-filling property}

Properties of space-filling curves are often different in finite and infinite approximations. \emph{Mazurkiewicz's theorem} states that there is no continuous, one-to-one mapping of a line segment onto a square, therefore no curves can be both self-avoiding and strictly space-filling, but a finite approximation of a given space-filling curve can be both self-avoiding and approximately space-filling. Our curves pass within a small distance from all points of the dark tiles which surround the curve and this distance can be arbitrarily reduced by carrying the recursive construction to an appropriate level. [PLF91]

\subsection{Edge-rewriting construction}

Here we show how to construct approximations of $F_n(k)$ fractal curves with ER method. Let us start with a generator curve $S_n$ on the overall grid which consists of $T_n$ edges on a triangular grid $(T_{n+1})$. All edges are lying on exactly one side of each dark tiles which stand in the same position as the big triangle ($F_n(0)=F_1$ generator pattern). Edges are originating from the left corner and terminating in the right corner of the big triangle. In the next approximation we divide all the dark tiles into $n$ equal pieces and we get $T_n^2$ dark subtriangles in the overall graph. We have to substitute the dark subtriangles with the contracted, rotated and reflected copies of the generator path in the direction of the edges. Direction codes define the tiles unambiguously. See \emph{Figure 4}.

\begin{remark}
For producing recursive curves we will use upper index $k$ at $w_n$ and $s_n$ direction code strings to mark the $k$-th approximation level.
\end{remark}

For example the simplest S-path means the following assignment: $s_2^0=0 \rightarrow s_2^1=105$ on the overall grid.

To any $p$ even and $q$ odd direction code the corresponding transformations (rotation and reflection) give the following 3 digits:

\begin{center}
$p \rightarrow \qquad (p+1) \mod 6, \qquad p, \qquad (p+5) \mod 6$

$q \rightarrow \qquad (q+5) \mod 6, \qquad q, \qquad (q+1) \mod 6$
\end{center}

The next approximation contains $T_n^2$ digits (9 edges). See \emph{Figure 8}.

\smallskip

\centerline{$s_2^1=105 \implies s_2^2=012105450$}

\begin{figure}[ht]
\centering
\includegraphics[width=0.8\textwidth]{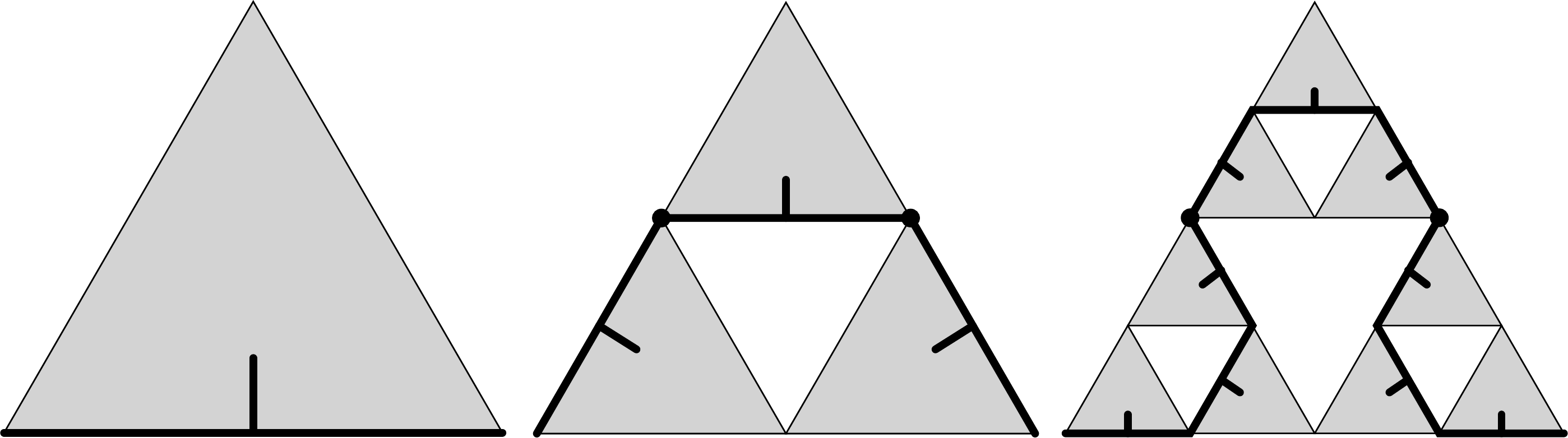}
\centerline{}
\centerline{\bf Figure 8. An S-path as an ER recursive curve.}
\centerline{\bf Sierpi\'{n}ski Arrowhead Curve $F_2(2)$.}
\centerline{\bf $0 \rightarrow 105 \implies 012105450$ ($s_2^0 \rightarrow s_2^1 \implies s_2^2$)}
\end{figure}

Let us see this assignment in general case from a direction code (1 digit as an even $p$ or an odd $q$ number) to $T_n$ edges (consisting of $5$ possible direction codes from set $x$ and set $y$).

The generalized formula for even numbers ($p$):

$p \rightarrow (p+x_1) \mod 6, (p+x_2) \mod 6, ... (p+x_{T_n}) \mod 6$, where $x_i \in \{0,1,4,5\}$

and the corresponding assignments for odd numbers ($q$) are:

$q \rightarrow (q+y_1) \mod 6, (q+y_2) \mod 6, ... (q+y_{T_n})\mod 6$ where $y_i \in \{0,5,2,1\}$.

Missing values in set $x$ and set $y$ mean 2 forbidden turns in S-path (no turn back to the same tile).

The order of the values in set $x$ and set $y$ indicates the corresponding pairs for the other parity, otherwise the reflected turns with the same angle from an even direction and from an odd direction. For example we have the following formula to an even direction: $(p+1) \mod 6$. The corresponding formula for an odd direction is $(q+5) \mod 6$ because the fourth value in set $x$ is a $1$ and in set $y$ it is a $5$.

\subsection{Node-rewriting construction}

We will use W-paths to get NR recursive curves. The simplest one is $w_2^0=0 \rightarrow w_2^1=15$ on the inscribed grid $(T_2)$, which has 3 nodes. In the next approximations we will keep all the previous edges, contracted them to the actual unit length and replace all the nodes among them to the contracted, rotated and reflected copies of our generator W-path as the last approximation of the S-path indicated how the dark tiles facing. See \emph{Figure 9}.

\begin{figure}[ht]
\centering
\includegraphics[width=0.8\textwidth]{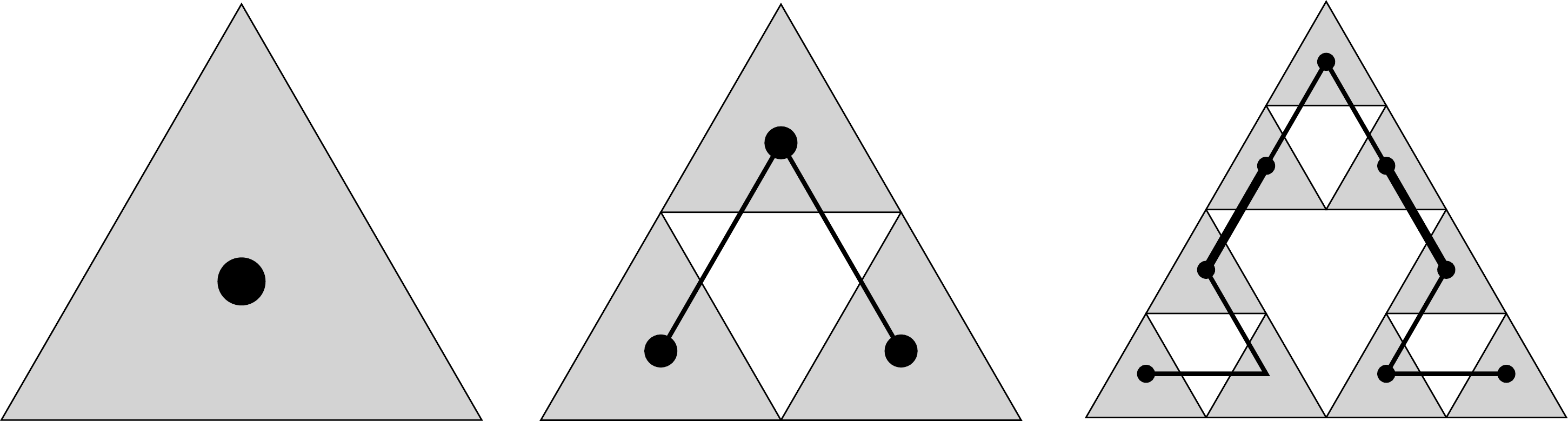}
\centerline{}
\centerline{\bf Figure 9. A W-path as an NR recursive curve.}
\centerline{\bf Sierpi\'{n}ski Arrowhead Curve $F_2(2)$.}
\centerline{\bf $0 \rightarrow 15 \implies 02115540$, ($w_2^0 \rightarrow w_2^1 \implies w_2^2$)}
\end{figure}

To get the next fractal approximation as a string $(w_n^{k+1})$ we have to insert these direction codes before, between and after the digits of our actual $w_n^k$ string but it is much more easy to transform the original W-path into an S-path first and using the ER method to construct the required $k$th approximation of the fractal and at the end, transform it back to NR code with \emph{Table 1}.

We get the following NR codes for the first 3 approximations (up to $T_2^3$ inscribed graph, $w_2^0 \rightarrow w_2^1 \implies w_2^2 \implies w_2^3$). Larger sized digits show the inheritance of the edges from the previous recursive levels:

\begin {center}
{$0 \rightarrow 15 \implies 02$ {\Large $1$} $15$ {\Large $5$} $40 \implies 15$ {\Large $0$} $02$ {\Large $2$}
$31$ {\Huge $1$} $02$ {\Large $1$} $15$ {\Large $5$} $40$ {\Huge $5$} $53$ {\Large $4$} $40$ {\Large $0$} $15$}
\end{center}

\subsection{Example for the least compound order}

We have shown how to construct recursive curves with both methods by W-paths and S-paths which are the generator curves of the NR and the ER recursive curves. 

\begin{figure}[ht]
\centering
\includegraphics[width=0.75\textwidth]{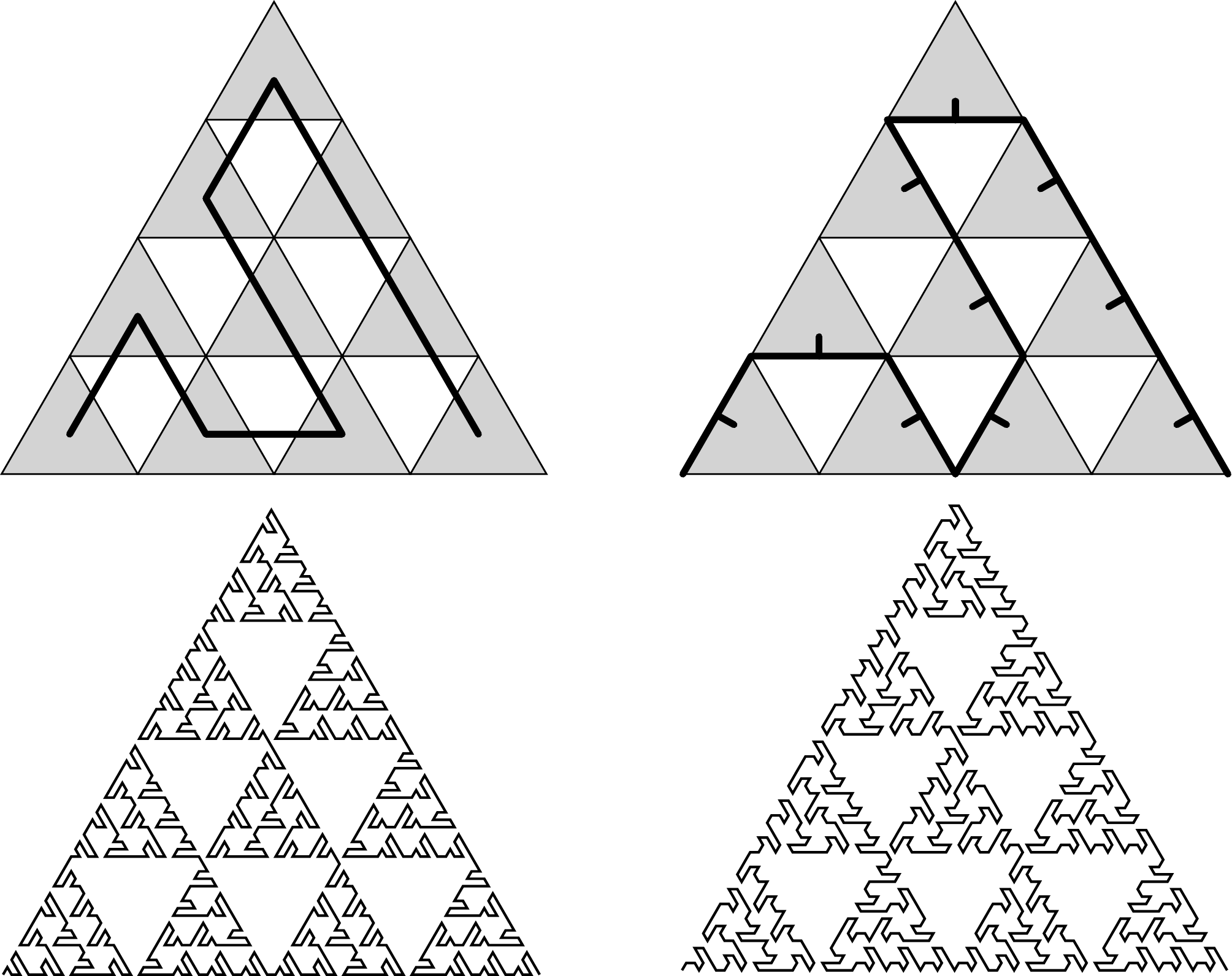}
\centerline{}
\centerline{\bf Figure 10. A bijective pair of $F_4$ generator curves (W- and S-paths)} 
\centerline{\bf and their $F_4(3)$ approximations (NR and ER recursive curves).}
\end{figure}

Both of our examples in \emph{Figures 8} and \emph{9} presented the first single result for $n=2$ which makes the unique bilateral symmetric permutation of the dark tiles in both methods and it keeps this property in all approximations. This curve is known as the \emph{Sierpi\'{n}ski Arrowhead Curve}. The commonly known ER version can be described easier but we can see the name of the curve came from the arrowhead shape of the NR generator curve. Every NR curve contains the original arrowhead shape at their tops.

This recursive curve family $F_n(\infty)$ for $n>2$ can be seen \emph{first time here}. See \emph{Figure 10} where you can compare the corresponding paths and recursive curves in both construction methods.

Corresponding recursive curves map the same fractal in the same permutation of the dark tiles in all approximations. This is an example from the first (least) compound order ($n=4$) where $F_n(\infty)$ fractal pattern first time differs from the \emph{Pascal Triangle modulo $n$} pattern.

See \emph{Table 2} for the numbers of the Hamiltonian-paths, which appear in [SEH05] and the numbers of the generator curves of the \emph{generalized Sierpi\'{n}ski Arrowhead Curve}, which \emph{first appear here as a new integer sequence}. We got these results by our computer program which is a smart backtrack algorithm.

The \emph{generalized Sierpi\'{n}ski Arrowhead Curve} fills the \emph{generalized Sierpi\'{n}ski Gasket} $(F_n(\infty))$ in both methods. Their \emph{Hausdorff dimension} [HD] is:

\smallskip

$D=log_n{(T_n)}=log_n{\dfrac{n(n+1)}{2}}$ \qquad where \qquad $\lim_{n\to\infty}log_n{(T_n)}=2$

\[
\begin{tabular}{|c||c||c||c|}
\hline
\bf
{\small n} & \bf {\small $T_n$} & \bf {\small $H_n$} & \bf {\small $W_n=S_n$ } \\ \hline \hline
\bf 2 & 3 & 1 &  1 \\ \hline
\bf 3 & 6 & 2 & 2 \\ \hline
\bf 4 & 10 & 10 & 4 \\ \hline
\bf 5 & 15 & 92 & 16 \\ \hline
\bf 6 & 21 & 1852 &  68 \\ \hline
\bf 7 & 28 & 78032 &  464  \\ \hline
\bf 8 & 36 & 6846876 &  3828 \\ \hline
\bf 9 & 45 & 1255156712 &  44488 \\ \hline
\end{tabular}
\]

\smallskip

\centerline{\bf \emph{Table 2.} Numbers of the Hamiltonian-paths $(H_n)$,}
\centerline{\bf the W-paths $(W_n)$ and the S-paths $(S_n)$}
\centerline{\bf on the generator pattern $F_n$ consisting of $T_n$ dark tiles.}

\section{Extendability of \emph{Theorem 1} to recursive curves}

In this section we will show why our bijection is also extendable to recursive curves.

Transformation between ER and NR recursive curves usually means their properties change. We will show our recursive curves unlike others always keep their self-avoiding property and all their edges keep the same length after transforming them into each other.

\subsection{Length of the connecting edges}

In [McK94] we can see tiling-paths of ER recursive curves which usually connect different facing tiles. In our case it means connections between different coloured tiles with tiling-paths (edges touch all used tiles only at one side). The connecting lengths are the same unit lengths but by transforming them to connected centroids we get two more connecting lengths. Centroids between edge-neighbour tiles are in $1/\sqrt{3}$ distance, between node-neighbour tiles they are in $2/\sqrt{3}$ distance from each other. Using S-path we have chosen appropriate connections. An S-path can connect only the same facing dark tiles which ensures the unit length of the connecting edges independently from the transformation between the paths.

\subsection{Self-avoiding property of our recursive curves}

ER curves based on S-paths. Consider the generator pattern as an $ABC$ triangle, where $A$ is the entering point, $B$ is the exiting point and $C$ is the top of the triangle. An S-path is always avoid corner $C$. By extending the \emph{overall grid} for larger approximations ($k>1$), the overall graph consists of $T_n^k$ little tile-like triangular grids. Each little triangle consists of $T_{n+1}$ grid points but they share common corners. Their edges always avoid $C$ corners of the little grids and they use the common corners between them only once as their exiting nodes ($B$) became to entering nodes ($A$) of different, neighbour tiles.

NR curves based on W-paths. We extend the inscribed grid to an inscribed graph which consists of $T_n^k$ little separated triangular grids with separated paths. We connect the little triangular grids with the edges of the W-path without the forbidden turns which ensures we stay inside this graph and we never use the corners of the little triangular grids twice. They could not be entering and also exiting points at the same time.

\subsection{Generator curves of the same fractal $F_n(\infty)$}

The bijective pairs of the W-paths and the S-paths as recursive curves fill and touch or pass over the dark tiles in all approximations exactly in the same order. They keep their self-avoiding property and continuously map the same fractal. We can transform them anytime into each other in the same way as the paths with \emph{Table 1}.

They approximate the same fractal $F_n(\infty)$ from inside (NR curve) and from outside (ER curve) and in an arbitrarily large finite $k$th approximation the distance of the paths becomes infinitely small anywhere. Contact points of the $ER$ curve and edges of the $NR$ curve are getting infinitely close to each other and we cannot distinguish the little triangles, their centroids and their touching edges any more.

By avoiding wrong turns in the construction of the recursive curves and by our other considerations described above the bijection is also valid for recursive curves.

\section {Transforming the paths into L-system codes}

We can transform our generator paths into the Lindenmayer-system to draw the fractal curves in an easy way. 

\subsection{Edge-rewriting codes}

We can transform our W-paths (NR generator curves) into ER strings in L-system by assigning the digit-pairs of the $w_n$ string with L-system symbols.

The first digit of the pair is an even $(p)$ or an odd $(q)$ number, then these assignings depend on the second digit of the pair:

\begin{multicols}{2}
$p \rightarrow A$ 

$(p+1) \rightarrow -B$

$(p+2) \mod{6} \rightarrow -B-$

$(p+5) \mod{6} \rightarrow A+$

$q \rightarrow B$

$(q+1) \mod{6} \rightarrow B-$

$(q+4) \mod{6} \rightarrow +A+$

$(q+5) \mod{6} \rightarrow +A$
\end{multicols}

For example on an $F_4$ generator pattern a possible W-path (supplemented with the leading and ending zero) is $w_4=01502215550$. Then by transforming the digit-pairs into (ER) L-system we get the following string: $A=-B+A+B--B-AA++A+BBB-$. Let the other assigning equal to its inverse: $B=+A-B-A++A+BB--B-AAA+$.

\[
\renewcommand{\arraystretch}{1.3}
\begin{tabular}{|c||c|c|c|c|c|c|}
\hline
\bf
$\blacksquare$ & {\bf $b=0$} & \bf {1} & \bf {2} & \bf {3} & \bf {4} & \bf {5} \\ \hline \hline
\bf {$a=0$} & A & -B & -B- & $\blacksquare$ & $\blacksquare$ & A+ \\ \hline
\bf 1 & +A & B & B- & $\blacksquare$ & $\blacksquare$ & +A+ \\ \hline
\bf 2 & $\blacksquare$ & A+ & A & -B & -B- & $\blacksquare$ \\ \hline
\bf 3 & $\blacksquare$ & +A+ & +A & B & B- & $\blacksquare$ \\ \hline
\bf 4 & -B- & $\blacksquare$ & $\blacksquare$ & A+ & A & -B \\ \hline
\bf 5 & B- & $\blacksquare$ & $\blacksquare$ & +A+ & +A & B \\ \hline
\end{tabular}
\]

\smallskip

\centerline {\bf \emph{Table 3.} Transformation table for changing the edge-pairs}
\centerline {\bf  of the W-paths to an ER string in L-system.}

\bigskip

You can see this W-path (NR generator curve) on the left side of \emph{Figure 10}. For the sake of simplicity we can use \emph{Table 3} to transform it to ER string in L-system. We have to read the table by rows and columns $(a,b)$. You can see this S-path (ER generator curve) on the right side of \emph{Figure 10}.

We can make various fractal patterns by changing the inverse rule (B=...) to other for example the palindrome of the B string or to the inverse of any other path in the same order of $n$.

\subsection{Node-rewriting codes}

To get the NR string in L-system we only have to insert F symbols among the symbol groups and we have to change all the A symbols to X, and all the B symbols to Y, because in L-system, A and B symbols are drawing instructions like F symbol (they draw the edges), but X and Y are only variables, we only need them to substitute the nodes which they symbolize.

Finally we get: $X=-YF+X+FY-F-Y-FXFX+F+X+FYFYFY-$ and its inverse $Y=+XF-Y-FX+F+X+FYFY-F-Y-FXFXFX+$ as strings in L-system which code the NR recursive curve.

For illustration see the left side of \emph{Figure 10}.

\subsection{Supplement}

L-system codes can be visualized easier with Inkscape [IS] freeware vector-graphic application program which has a built-in L-system extension, or in an online L-system application [OL]. In this paper we discussed recursive curves with simple axioms against for example Moore-curve, but it is also possible to make cycles from them or combine them with compound axioms.

For further information and experiments with recursive curves we recommend that beyond the mentioned sources curious readers read the beautiful book of Jeffrey Ventrella [V12] and [A16]. Gosper-curve variants are also related [FSN01].

\section* {Summary}

We have observed the properties of Hamiltonian-paths and tiling-paths on the triangular grid and their transformability into each other. We have found and proved the bijective relation of their subsets in \emph{Section 1}.

In \emph{Sections 2} and \emph{3} we have used our bijective related assymetric paths as generator curves by which we have constructed the \emph{generalized Sierpi\'{n}ski Arrowhead Curve} in larger graphs also in node-rewriting and edge-rewriting method. The cardinality of their generator curves specifies a new integer sequence. 

In \emph{Section 4} we have presented a transformation table to change the direction code strings of the paths into L-system strings to draw the fractals.

Our previous papers can complete this field with some details. For all of my papers please check my Google Scholar site [KA].

\section* {Acknowledgement}
We would like to express our thanks to \emph{Mih\'aly Hujter} professor of Mathematics at TU Budapest, to \emph{Guszt\'av Ga\'al} mathematician at E\"{o}tv\"{o}s Lor\'and University (ELTE) Budapest and the inspiring community of the \emph{Mathematics Museum in Budapest}.

\section* {References}

\smallskip \noindent [McK94] McKenna, Douglas M.: \emph{SquaRecurves, E-Tours, Eddies and Frenzies: Basic Families of Peano Curves on the Square Grid}, In: Guy, Richard K., Woodrow, Robert E.: \emph{The Lighter Side of Mathematics: Proceedings of the Eugene Strens Memorial Conference on Recreational Mathematics and its History}, pp. 49-73, Mathematical Association of America, 1994.
\smallskip

\smallskip \noindent [PL90] Prusinkiewicz, P. and Lindenmayer, A.: \emph{The algorithmic beauty of plants}, Springer, 1990.
\smallskip 

\smallskip \noindent [PLF91] Prusinkiewicz, P., Lindenmayer, A., and Fracchia, F. D.: \emph{Synthesis of space---filling curves on the square grid}, In: H.-O. Peitgen, J. M. Henriques and L. F. Penedo, eds.: \emph{Fractals in the fundamental and applied sciences}, pp. 341-366, North-Holland, 1991. 
\smallskip 

\smallskip \noindent [Sa94] Sagan, H.: \emph{Space---filling curves}, Springer, 1994.
\smallskip 

\smallskip \noindent [M82] Mandelbrot, B. B.: \emph{The Fractal Geometry of Nature}, W. H. Freeman and Company, New York, 1982.
\smallskip 

\smallskip \noindent [SEH05] Staji\'c, J., Elezovi\'c-Had\v{z}i\'c, S.: \emph{Hamiltonian walks on Sierpinski and n-simplex fractals}, \url{https://arxiv.org/abs/cond-mat/0310777}, 2005.
\smallskip 

\medskip \noindent [HD] \emph{List of fractals by Hausdorff dimension}, \url{https://en.wikipedia.org/wiki/List_of_fractals_by_Hausdorff_dimension}
\smallskip

\smallskip \noindent [IS] \emph{Inkscape vector graphic application with built in L-system}, \\
\url{https://inkscape.org/en/}
\smallskip

\smallskip \noindent [OL] \emph{Online Lindenmayer-system application}, \\
\url{http://www.kevs3d.co.uk/dev/lsystems/}
\smallskip

\smallskip \noindent [V12] Ventrella, J.: \emph{Brainfilling Curves --- a Fractal Bestiary}, \\
\url{http://www.brainfillingcurves.com/}, 2012.
\smallskip

\smallskip \noindent [A16] Arndt, J.: \emph{Plain--filling curves on all uniform grids}, \\
\url{https://arxiv.org/pdf/1607.02433.pdf}, 2016.
\smallskip

\smallskip \noindent [FSN01] Fukuda, H., Shimizu, M., Nakamura, G.: \emph{New Gosper Space Filling Curves}, In: Proceedings of the International Conference on Computer Graphics and Imaging (CGIM2001) \\
\url{http://kilin.clas.kitasato-u.ac.jp/museum/gosperex/343-024.pdf}, 2001.

\smallskip \noindent \rm [KA] \emph{Google Scholar citations of Kaszanyitzky, A.},
\\ \url{https://scholar.google.hu/citations?user=i5daxSoAAAAJ}
\smallskip

\end{document}